\documentclass[12pt]{amsart}

\usepackage{lmodern}
\usepackage[a4paper,left=2.5cm,right=2.5cm,top=2.5cm,bottom=2.5cm]{geometry}
\usepackage[T1]{fontenc}
\usepackage{amsmath,amssymb,mathtools}
\usepackage[colorlinks=true,linkcolor=blue,citecolor=blue,urlcolor=blue]{hyperref}

\newcommand{\calG}{\mathcal G}
\newcommand{\Z}{\mathbb Z}

\title[A Note on the 2-Local Homotopy Types of $G_2$-Gauge Groups]{A Note on the 2-Local Homotopy Types\\ of $G_2$-Gauge Groups}

\author{Phuc Vo Dang}
\address{Department of Mathematics, FPT University, An Phu Thinh New Urban Area, Quy Nhon, Vietnam}
\email{dangphuc150488@gmail.com}
\thanks{ORCID: \url{https://orcid.org/0000-0002-6885-3996}}


\begin{document}

\begin{abstract}
In a recent preprint \cite{Kameko2026} , Kameko presented a substantial completion of the 2-local classification of $G_2$-gauge groups over $S^4$, extending earlier work by Kishimoto, Theriault, and Tsutaya. The central strategy relies on reducing the 2-local classification to the order of the Samelson product $\langle i_3,1\rangle$ and separating specific gauge group homotopy types. The purpose of this note is to provide necessary mathematical refinements and localization clarifications to several key steps in the proof of the classification theorem. Specifically, we refine the integral isomorphism claim for gauge group homotopy to its correct 2-local form, resolve an EHP sequence extension regarding the injectivity of the Hopf invariant, and make explicit the Postnikov layer conventions required for the mod 2 Leray--Serre spectral sequence calculations. We confirm that with these adjustments, the main 2-local classification theorem in \cite{Kameko2026} also holds as claimed.
\end{abstract}

\maketitle

\section{Introduction}

The classification of gauge groups over spheres is a central problem in unstable homotopy theory. In a recent preprint \cite{Kameko2026}, Kameko gave a convincing completion of the 2-local part of the classification of $G_2$-gauge groups over $S^4$, continuing the foundational work of Kishimoto, Theriault, and Tsutaya \cite{KishimotoTheriaultTsutaya2017}. The overarching strategy (namely, the reduction of the 2-local classification to the order of the Samelson product $\langle i_3,1\rangle$ and the separation of $\calG_4$ from $\calG_0$) is natural and compatible with Lang's identification of the boundary map with the triple adjoint of a Samelson product \cite{Lang1973}. Furthermore, the computations utilizing Toda's notation for $\bar\nu_6$, $\epsilon_6$, $\zeta_6$, the relative Samelson product calculations based on James' formula \cite{James1976}, and the parity arguments involving $p_{11}q$ are structurally robust, drawing upon classical homotopy-group inputs from Mimura \cite{Mimura1967} and Toda \cite{Toda1962}.

The objective of this note is to strengthen the theoretical foundation of the classification in \cite{Kameko2026} by addressing two mathematical inaccuracies in the intermediate arguments, providing a necessary localization clarification regarding connective covers, and refining several technical expressions. It should be emphasized that these modifications do not invalidate the main conclusion of \cite{Kameko2026}; rather, they place the proof of the main theorem on completely rigorous ground.

\section{Integral vs. 2-Local Homotopy of Gauge Groups}

In \cite[Proposition 6.2]{Kameko2026}, it is asserted integrally that the induced homomorphism
\[
 h_{k*}\colon \pi_3(\calG_k)\longrightarrow \pi_3(G_2)
\]
is an isomorphism. The argument provided relies on the vanishing of the homotopy group $\pi_3(\Omega^3_0G_2) \cong \pi_6(G_2)$. However, as established in Mimura's classical calculations \cite{Mimura1967} (and recorded in Section 2 of \cite{Kameko2026}), one has integrally
\[
 \pi_6(G_2)\cong \Z/3.
\]
Consequently, the integral long exact sequence associated with the evaluation fiber sequence yields
\[
 0\longrightarrow \pi_3(\calG_k)\xrightarrow{h_{k*}}\pi_3(G_2)\longrightarrow \pi_6(G_2),
\]
and the integral isomorphism conclusion does not follow directly without further analysis of the boundary map $\pi_3(G_2)\to \pi_6(G_2)$.

To establish the classification, the integral statement is not strictly necessary; rather, one requires the corresponding 2-local assertion. We formulate the correct statement as follows:

\begin{description}
\item[Claim] The 2-localized homomorphism
\[
 h_{k*}\colon \pi_3(\calG_k)_{(2)}\longrightarrow \pi_3(G_2)_{(2)}
\]
is an isomorphism.
\end{description}

Indeed, applying 2-localization to the long exact homotopy sequence of the fibration
\[
 \calG_k\xrightarrow{h_k}G_2\xrightarrow{\partial_k}\Omega^3_0G_2,
\]
we obtain the exact sequence
\[
 \pi_7(G_2)_{(2)}\longrightarrow \pi_3(\calG_k)_{(2)}
 \xrightarrow{h_{k*}}\pi_3(G_2)_{(2)}
 \longrightarrow \pi_6(G_2)_{(2)}.
\]
Since $\pi_7(G_2)=0$ and $\pi_6(G_2)=\Z/3$, both outer 2-local groups vanish identically. Exactness immediately guarantees the required 2-local isomorphism. This refined 2-local form is precisely what is utilized in the subsequent arguments of \cite{Kameko2026}, specifically when the proof of Proposition 1.2 applies the result to
\[
 (h_4\circ\varphi)_*\colon \pi_3(G_{2,11})_{(2)}\longrightarrow \pi_3(G_2)_{(2)}.
\]

\section{Injectivity of the Hopf Invariant and EHP Sequences}

A second refinement is required in the first paragraph of the proof of Proposition 2.2 in \cite{Kameko2026}, where the vanishing of $\pi_{10}(S^5)$ is assumed in order to infer the injectivity of the Hopf invariant
\[
 H\colon \pi_{11}(S^6)_{(2)}\longrightarrow \pi_{11}(S^{11})_{(2)}.
\]
According to Toda's tables \cite{Toda1962}, the integral homotopy group is non-zero:
\[
 \pi_{10}(S^5)\cong \Z/2.
\]
While this observation affects the immediate rationale, the intended injectivity conclusion remains valid and can be rigorously established through the following modified argument. 

Consider the relevant portion of the EHP exact sequence:
\[
 \pi_{10}(S^5)_{(2)}\xrightarrow{\Sigma}\pi_{11}(S^6)_{(2)}\xrightarrow{H}\pi_{11}(S^{11})_{(2)}.
\]
Here, $\pi_{10}(S^5)_{(2)}\cong \Z/2$, whereas $\pi_{11}(S^6)_{(2)}\cong \Z_{(2)}\{\Delta(\iota_{13})\}$ is torsion-free. Since any homomorphism from a torsion group ($\Z/2$) to a torsion-free group ($\Z_{(2)}$) must be the zero map, the suspension homomorphism $\Sigma$ vanishes identically. Exactness therefore implies $\ker H=0$, proving that $H$ is indeed injective on $\pi_{11}(S^6)_{(2)}$. 

Because the localization map $\pi_{11}(S^6)\cong\Z\to \pi_{11}(S^6)_{(2)}$ is injective, this 2-local injectivity is entirely sufficient to deduce the subsequent integral equality in the infinite cyclic group $\pi_{11}(S^6)$. The Hopf invariant computation establishing
\[
 (2\iota_6)\circ \Delta(\iota_{13})\simeq 4\Delta(\iota_{13})
\]
then proceeds without alteration, preserving its application in Proposition 3.9 of \cite{Kameko2026}.

\section{Postnikov Layers and Connective Covers}

In Section 5 of \cite{Kameko2026}, a clarification is necessary regarding the sequence of connective covers
\[
 G_2\langle 9\rangle\xrightarrow{g_3}G_2\langle 8\rangle\xrightarrow{g_2}G_2\langle 3\rangle\xrightarrow{g_1}G_2,
\]
where the successive fibers are stated to be $K(\Z,2)$, $K(\Z/2,7)$, and $K(\Z/2,8)$. Integrally, this description is incomplete: the passage from the 3-connected cover to the 8-connected cover must also eliminate the 3-primary group $\pi_6(G_2)\cong \Z/3$, and the passage from the 8-connected cover to the 9-connected cover detects the odd-primary component of $\pi_9(G_2)\cong \Z/6$. Consequently, the stated Eilenberg--Mac Lane fibers are accurate only after restricting to the 2-primary Postnikov layers, or working modulo the Serre class of odd torsion groups.

This distinction does not impact the mod 2 spectral sequence calculations, as Eilenberg--Mac Lane spaces attached to odd-primary groups exhibit trivial reduced mod 2 homology. To make the proof of Proposition 3.10 completely rigorous, one must explicitly adopt the convention at the beginning of Section 5 that the relevant Postnikov layers are considered after 2-localization---equivalently, that the Leray--Serre spectral sequence is evaluated after discarding odd-primary fibers, which are mod 2 acyclic. Under this convention, the transgressions
\[
 d_8(u_7)=y_8,
 \qquad d_9(Sq^1u_7)=Sq^1y_8,
 \qquad d_{11}(Sq^2Sq^1u_7)=Sq^2Sq^1y_8
\]
are precisely the relevant ones up to total degree $12$, and the deduction
\[
 (g_1\circ g_2\circ g_3)_*\bigl(H_{11}(G_2\langle 9\rangle;\Z)/\operatorname{Torsion}\bigr)
 \subset 8\,H_{11}(G_2;\Z)
\]
remains fully justified, supplying the necessary divisibility input for Proposition 3.10.

\section{Further Refinements and Notational Corrections}

For the completeness of the mathematical record, we note several smaller technical and notational adjustments that refine the arguments in \cite{Kameko2026}:

\subsection*{Cofiber Sequences in Proposition 2.4}
In \cite[Proposition 2.4]{Kameko2026} , the assertion $$[P^{13}(G_2),G_2]\cong 0$$ should be corrected to read
\[
 [P^{13}(2),G_2]\cong 0.
\]
This notation aligns with the preceding cofiber sequence $S^{12}\to P^{13}(2)\to S^{13}$.

\subsection*{Parity and Multiplication Factors in Proposition 2.5}
In  \cite[Proposition 2.5]{Kameko2026}, the conclusion drawn from the parity of $a$ should not be written as $(\psi\circ s'')\simeq 4\langle \bar\nu_6+\epsilon_6\rangle$. Because the calculation evaluates $84(\psi\circ s'')$, the conclusion is more accurately expressed as
\[
 84\cdot(\psi\circ s'')\simeq 4\cdot\langle \bar\nu_6+
 \epsilon_6\rangle
\]
after suppressing the odd unit in the cyclic 2-primary summand. More explicitly, the proof yields
\[
 84\cdot(\psi\circ s'')\simeq 4a\cdot\langle \bar\nu_6+\epsilon_6\rangle
\]
with $a$ odd, which is the exact form required upon precomposition with $\Sigma^3(p_{11}q)$.

\subsection*{Typographical Corrections in Propositions 3.3 and 4.3}
In  \cite[Proposition 3.3]{Kameko2026}, the displayed expression $\langle\bar\nu_6+\epsilon\rangle$ should explicitly read $\langle\bar\nu_6+\epsilon_6\rangle$. Similarly, in  \cite[Proposition 4.3]{Kameko2026}, the final displayed computation for the case $\gamma=\epsilon_6$ should begin with $\langle i_3,\epsilon_6\rangle_r$ rather than $\langle i_3,\bar\nu_6\rangle_r$. The stated result $\langle i_3,\epsilon_6\rangle_r\simeq 2\bar\nu_6\circ\nu_{14}$ is unaffected.

\subsection*{Sharpening Null-Homotopy Statements in Proposition 3.6}
In  \cite[Proposition 3.6]{Kameko2026}, the phrase ``\textit{By Proposition 3.1, it is null homotopic}'' warrants sharpening. Proposition 3.1 implies that the class of
\[
 21\,p\psi\Sigma^3(q\langle\bar\nu_6+\epsilon_6\rangle)
\]
vanishes after quotienting $2\pi_{17}(S^6)$ by the direct summand generated by $\bar\nu_6\circ\nu_{14}$. This quotient step is essential to force the $\zeta_6$-coefficient to vanish modulo $8$. The element itself is not zero prior to passing to this quotient, as $\bar\nu_6\circ\nu_{14}$ has order $4$.

\subsection*{Retention of Numerical Factors in Proposition 1.5 and Section 6}
In the proof of \cite[Proposition 1.5]{Kameko2026}, the transition from Propositions 3.3 and 3.9(2) to the displayed identity
\[
 p\circ\psi\circ s''\circ\Sigma^3(p_{11}\circ q\circ\langle\bar\nu_6+
 \epsilon_6\rangle)
 \simeq a''(\bar\nu_6+\epsilon_6)\circ c\nu_{14}
\]
should explicitly account for the factor $21$ from Proposition 3.3. While this omission is mathematically harmless because the target summand generated by $\bar\nu_6\circ\nu_{14}$ has order $4$ and $21\equiv 1\pmod 4$, retaining the factor ensures formal precision. 

Similarly, in \cite[Proposition 6.4]{Kameko2026}, the reference back to Section 3 should include this factor: the result established in Section 3 is that there exists an odd integer $a''$ such that
\[
 21\cdot(p\circ\psi\circ s'')\simeq a''(\bar\nu_6+\epsilon_6),
\]
rather than $p\circ\psi\circ s''\simeq a''(\bar\nu_6+\epsilon_6)$. With this correction, the displayed identity
\[
 21\cdot 3\cdot \partial_{4*}(i')\simeq 84\cdot\psi\circ s''\circ\Sigma^3(p_{11}\circ q')
\]
connects seamlessly to the argument of Proposition 2.5, confirming that $21\cdot 3\cdot\partial_{4*}(i')$ has order $2$.

\subsection*{Type-Correctness in Proposition 6.1}
In \cite[Proposition 6.1]{Kameko2026}, the displayed target of the homomorphism $\partial_{1*}$ should be $\pi_{17}(G_2)$, or equivalently $[\Sigma^3S^{14},G_2]$, rather than $[\Sigma^3G_{2,11},G_2]$. Since the elements to which $\partial_{1*}$ is applied in that paragraph belong to $\pi_{14}(G_2)$, writing the target as $[\Sigma^{3}S^{14}, G_{2}]$ preserves type-correctness.

\subsection*{Grammatical and K-theory Notations in Proposition 5.2}
In \cite[Proposition 5.2]{Kameko2026}, the formulation should explicitly specify ``\textit{the mod 2 reduced cohomology group of $G_2\langle 8\rangle$ is isomorphic to}'' and ``\textit{up to degree $\leq 12$}''. Furthermore, the final sentence should reference $H^{11}(G_2\langle 8\rangle;\Z)$ instead of $H^{11}(G\langle 8\rangle;\Z)$.

\section{Conclusion}

By incorporating the 2-local refinement of \cite[Proposition 6.2]{Kameko2026}, establishing the injectivity of the Hopf invariant via the torsion-free property of $\pi_{11}(S^6)_{(2)}$, and adopting explicit 2-local Postnikov conventions in Section 5, the logical architecture of \cite{Kameko2026} becomes entirely rigorous. \cite[Proposition 6.1]{Kameko2026} successfully distinguishes $\calG_1$ and $\calG_2$ from $\calG_0$ and $\calG_4$ through the computation
\[
 \pi_{13}(\calG_k)\cong \Z/(8,2k)\oplus\Z/2.
\]
The separation of $\calG_4$ from $\calG_0$ relies strictly on 2-local information on $\pi_3$, for which the 2-localized Proposition 6.2 is completely adequate. Consequently, the determination of the 2-primary order of $\langle i_3,1\rangle$ stands, and the complete 2-local classification by $(k,8)$ follows as claimed in \cite[Theorem 1.1]{Kameko2026}.

\end{document}